# Leadership exponent in the pursuit problem for 1-D random particles

## G. Molchan

Institute of Earthquake Prediction Theory and Mathematical Geophysics,

Russian Academy of Science, 84/32 Profsoyuznaya st.,

117997, Moscow, Russian Federation

E-mail address: molchan@mitp.ru

**Abstract**. For n + 1 particles moving independently on a straight line, we study the question of how long the leading position of one of them can last. Our focus is the asymptotics of the probability p(T,n) that the leader time will exceed T when n and T are large. It is assumed that the dynamics of particles are described by independent, either stationary or self-similar, Gaussian processes, not necessarily identically distributed. Roughly, the result for particles with stationary dynamics of unit variance is as follows: L= -log p(T,n) /(Tlog n)=1/d+o(1), where $d/(2\pi)$ is the power of the zero frequency in the spectrum of the leading particle, and this value is the largest in the spectrum. Previously, in some particular models, the asymptotics of L was understood as a sequential limit first over T and then over n. For processes that do not necessarily have non-negative correlations, the limit over T may not exist. To overcome this difficulty, the growing parameters T and n are considered in the domain clog T<log n<CT, where c>1 . The Lamperti transform allows us to transfer the described result to self-similar processes with the normalizer of log p(T,n) becoming log T Log n.

## 1. Introduction and the main result.

In this paper we consider the pursuit problem in an ensemble of particles with random dynamics. The problem involves a population of particles on a straight line consisting of a single "pursued" particle and $n$ "pursuing" particles. The motion of the particles is described by independent random Gaussian processes. It is assumed the pursued particle is ahead of the others at the start. The main problem is the distribution of the time $\tau_n$ that it takes to catch up with the pursued particle.

The problem was the subject of a lively discussion for particles with uniform Brownian dynamics in the 1990-2010s. The issue as to the finiteness of the mean $E\tau_n$ has turned out to be a nontrivial problem. For the Brownian particles De Blassie [6,7] (see also [3]) has shown that

$$P(\tau_n > T) \sim cT^{-\gamma_n}, T \to \infty \tag{1.1}$$

where $c$ depends on the initial particle positions,

$$2\gamma_n = \sqrt{\lambda_1 + (n-1)^2/4} - (n-1)/2$$

and $\lambda_1$ is the first (principal) eigenvalue of the Dirichlet problem (with zero boundary conditions) for the Laplace-Beltrami operator on the subset $G_{n+1} \cap S^n$ of the unit sphere $S^n \subset R^{n+1}$:

$$G_{n+1} = \{x = (x_0,..., x_n) : x_i - x_0 < 0, i = 1,..., n\}.$$

The estimation of $\lambda_1$ is a technically complicated problem; nevertheless, it has been shown working on these lines that the mean $E\tau_n$ is infinite when $n \leq 3$ [4] and that it is finite when $n = 4$ [20]. Earlier, the finiteness of $E\tau_n$ for $n \geq 5$ has been proved by another method, [13].

Kesten [10] raised the issue of the asymptotics of $\gamma_n$ and showed that in the Brownian case

$$P(\tau_n > T) > T^{-\ln n \cdot (1+\varepsilon)/4}, T > T_0(n), n > n_0(\varepsilon) \tag{1.2}$$

for any $\varepsilon > 0$.

For diffusing particles, Krapivsky and Redner [11] made strong arguments in favor of the asymptotics $\gamma_n = \ln(4Rn)/(4R)$, where $R = \sigma_0/\sigma$ is the ratio of the diffusivity of the leading particle, $\sigma_0$, to the corresponding parameter for the other particles, $\sigma$.

Li and Shao [14, 15] considered a non-Markov model in which the particle dynamics is described by the fractional Brownian motion (FBM) $B_H(t)$, i.e., by the Gaussian process with stationary increments of the form

$$E(B_H(t) - B_H(s))^2 = |t-s|^{2H}, 0 < H < 1, B_H(0) = 0. \tag{1.3}$$

The case $H = 1/2$ corresponds to Brownian motion.

The Li and Shao result is as follows: for independent processes $B_H^{(i)}(t), i \geq 0$,

$$P\{B_H^{(i)}(t) - B_H^{(0)}(t) < 1, 0 < t < T, 1 \leq i \leq n\} = T^{-\gamma_{n,H}+o(1)}, \tag{1.4}$$

as $T \to \infty$ and

$$1/d_H \le \liminf_{n \to \infty} \gamma_{n,H} / \ln n \le \limsup_{n \to \infty} \gamma_{n,H} / \ln n < \infty \ ,$$

where

$$d_H = 2\int_0^\infty [e^{tH} + e^{-tH} - (e^{t/2} - e^{-t/2})^{2H}]dt = \frac{2\Gamma(1-H)\Gamma(2H)}{\Gamma(1+H)} \ . \tag{1.5}$$

The two-sided estimates of the exponent for the Brownian particles coincided, namely $\gamma_{n,1/2} = \ln n /(4 + o(1))$. This enabled Li and Shao to conjecture that for the general case of $H$ we must have

$$\gamma_{n,H} = \ln n /(d_H + o(1)) \ . \tag{1.6}$$

Our goal is to consider the case of particles with heterogeneous dynamics, and to use the limit transition in which T and n are varied simultaneously. In this way, we can prove the asymptotics (1.4, 1.6).

More specifically, we will consider an ensemble of independent continuous centered Gaussian processes $\{X^{(i)}(t), i \ge 0\}$. This may be a family of either stationary or self-similar processes. An example of the latter is the fractional Brownian motion. Our problem is the asymptotics of the probability

$$p_{T,n}(c) = P\{X^{(i)}(t) - X^{(0)}(t) < c, i = 1,..., n, 0 \le t \le T\}$$

as $T, n \to \infty$. We will consider $c = 1$ in the conditional situation $\{X^{(i)}(0) = 0, i \ge 0\}$ and $c = 0$ in the unconditional one.

The expected order of the $p_{T,n}(0)$ asymptotics in the stationary case can be illustrated by the following simple example. Consider $(n+1)$ i.i.d. Gaussian sequences $X^{(i)}(j), 1 \le j \le S$ with i.i.d. components. Then

$$p_{S,n}(0) = P\{X^{(i)}(j) < X^{(0)}(j), i = 1,..., n, 0 \le j \le S\} = (n+1)^{-S} \approx e^{-S \ln n}$$

In the general stationary case, it is natural to expect that the $p_{T,n}(0)$ for large $T, n$ depends mainly on the described array of random variables hidden in the data with $S = T/d$. To implement this idea, we relied on the works [2,8,14] and on the theory of extreme values for random Gaussian processes, [16].

**Theorem 1**. Let $\{X^{(i)}(t), i \ge 0\}$ be the ensemble of Gaussian stationary processes with correlation functions $\{r_i(t), r_i(0) = 1\} \subset R$. Assume that the $R$ set is finite, $r_0(t) \in L_1$ and $d_0 = \int r_0(t)dt > 0$ (in particular, $X_0(t)$ has a continuous spectral density $f_0(\lambda)$).

***The lower bound***. Let a) $1 - r_0(t) \le c|t|^{2H_0}, 0 < H_0 < 1 . |t| < \varepsilon_0$ Then

$$\liminf_{T,n \to \infty} \ln p_{T,n}(0) /(T \ln n) \ge -1/d_0 \tag{1.7}$$

provided that $T, n$ increase so that $\ln n \leq CT$.

**The upper bound.** Let for any $i \geq 1$

b) $1 - |r_i(t)| \geq c_i |t|^{2h_i} \wedge \delta_i, 0 < \delta_i < 1, h_i > 0, c_i > 0$;

c) $|r_i(t)| \leq \tilde{r}_i(t) \in L_1$, where $\tilde{r}_i(t)$ is a monotone function;

d) $f_0(\lambda) \leq f_0(0) < \infty$, e) $f_0(\lambda) \leq \varphi(\lambda) \in L_1$, where $\varphi(\lambda)$ is a monotone function,

and f) $\ln f_0(\lambda) \in L_{1,loc}$.

Then

$$\limsup_{T,n \to \infty} \ln p_{T,n}(0)/(T \ln n) \leq -1/d_0, \qquad (1.8)$$

provided that $T, n$ increase so that $c \ln T \leq \ln n$, with any fixed $c > 1$.

**Remark 1**. Theorem 1 allows us to consider the value $1/d_0$ as a *leadership exponent*. The following factors are crucial for the result (1.7, 1.8): the dynamics of the leading particle, the independence of the remaining particles, and the variance equality: $r_i(0) = const$, $i \geq 0$.

**Remark 2**. Obviously, $d_0 = \int r_0(t)dt = 2\pi f_0(0)$ and $2\pi f_0(\lambda) \leq \int |r_0(t)|dt$. Therefore the (d)-condition ($f_0(\lambda) \leq f_0(0)$) is automatically met if the correlation function of the leading particle is non-negative, $r_0(t) \geq 0$. In addition, in this case a modified normal comparison inequality proposed in [14] allows us to prove the upper bound (1.8) without any restrictions on the growth of the parameters $T$ and $n$. The proof is omitted because it does not require new ideas but complicates the text.

Finiteness of the $R$ set is not essential. This condition is used only in order to have a uniform restrictions of the parameters in the stochastic dynamics of particles.

**Remark 3**. The above pursuit problem can be interpreted as a *persistence* problem for the random field $\xi(t,i) = X^{(i)}(t) - X^{(0)}(t)$ on the expanding region $[0,T] \times [1,n]$. Only a few nontrivial examples of Gaussian fields are known for which the persistence exponent (in our case $1/d_0$) is obtained explicitly (see e.g. [19]). In these examples the expanding is controlled by a one-dimensional parameter, say, based on the similarity of the regions. In our case $\log p_{T,n}(o) = O(T \log n)$. Therefore the restriction of the type $c < \log n/T < C$ for the stationary case can by considered as an analogous of the similarity of the regions. This restriction is included in the conditions of Theorem 1. The background of the persistence problem for random processes can be found in the surveys [1, 5].

*Self-similar processes*. Originally, the pursuit problem was considered for (fractional) Brownian particles whose trajectories are independent copies of a special H-self-similar (H-ss) process. By definition, for any H-ss process $\{X(t), X(0)=0\}$, all finite-dimensional distributions of $X(t)$ and $\lambda^{-H} X(\lambda t)$ are identical for any $\lambda > 0$, shortly $X(t) \stackrel{d}{=} \lambda^{-H} X(\lambda t)$. The Lamperti transformation, $\tilde{X}(\tau) = X(e^\tau) e^{-\tau \cdot H}$, maps the H-ss process $\{X(t), 1 \leq t \leq T\}$ to a stationary process $\{\tilde{X}(\tau), 0 \leq \tau \leq \ln T\}$. As a result, for any H-ss processes $\{X^{(i)}(t)\}$ we shall have

$$P(M_n([1,T]) \leq 0) = P(\tilde{M}_n([0, \ln T]) \leq 0) = \tilde{p}_{\ln T, n}(0),$$

where $M_n(\Delta) = \max\{X^{(i)}(t) - X^{(0)}(t), t \in \Delta, 1 \leq i \leq n\}$, and the symbol "~" is related to $\tilde{X}^{(i)}(\tau)$.

The pursuit problem for the H-ss processes on $[0,T]$ is considered as an asymptotics of probability $p_{T,n}(1) = P(M_n([0,T]) \leq 1)$. Theorem 1 can be useful in this case if the log-asymptotics of $P(M_n([0,T]) \leq 1)$ and $P(M_n([1,T]) \leq 0)$ are identical. The next statement shows that these asymptotics with the normalizing factor $[\log T \log n]^{-1}$ exist only simultaneously and are identical. This kind of statement is quite expected (see, for example, [18]), but it needs proof to have constraints on the growth of $(T,n)$-parameters.

**Theorem 2. 1)** Let $\{X^{(i)}(t), i \geq 0, t \geq 0\}$ be a family of the H-ss processes with finite set of different correlation functions $\{r_i(t,s)\}$, $\psi_{T,n} = \log T \log n$.

1. Assume that $T, n(T) \to \infty$ in such a way that $\ln n \leq C \ln T$, then

$$\liminf \ln P(M_n([1,T]) \leq 0)/\psi_{T,n} \leq \liminf \ln P(M_n([0,T]) \leq 1)/\psi_{T,n}. \quad (1.11)$$

2. Assume that

$$r_0(t,1) t^{-H-1} \in L_1(0, \infty) \text{ and } d_0 = \int_0^\infty r_0(t,1) t^{-H-1} dt > 0. \quad (1.12)$$

Then

$$\limsup \ln P(M_n([1,T]) \leq 0)/\psi_{T,n} \geq \limsup \ln P(M_n([0,T]) \leq 1)/\psi_{T,n}, \quad (1.13)$$

where both limits are related to an arbitrary but common sequence $T, n(T) \to \infty$.

**Corollary 3.** The pursuit problem for fractional Brownian particles with $X^{(i)}(t) \stackrel{d}{=} B_H(t), i = 0, ..., n$ has the asymptotics

$$\log p_{T,n}(1) /[\log T \log n] = -1/d_H (1 + o(1)) \quad (1.14)$$

with the hypothetical Li-Shao constant $d_H$, (1.5), provided that $c \ln \ln T < \ln n \leq C \ln T, c > 1$ and $T \to \infty$. In addition, the result remains valid after replacing $\{B_H^{(i)}(t), i=1,\ldots,n\}$ with independent processes $\{t^{H-H_i}B_{H_i}^{(i)}(t), i=1,\ldots,n\}$, where $0 < H, \min H_i, \max H_i < 1$ and the number of different $H_i$ is finite.

***Remark 4.*** Following Remark 2, we note that $-1/d_H$ as the upper bound of the left part (1.14) can be obtained without additional restrictions on the growth of the $(T,n)$ parameters.

***Proof of Corollary 3.*** The Lamperti transform of the fractional Brownian motion $B_H(t)$ is the stationary process $\tilde{B}_H(\tau)$ with the correlation function:

$$r_H(t) = [e^{tH} + e^{-tH} - (e^{t/2} - e^{-t/2})^{2H}]/2 > 0.$$

Now the conditions of Theorem 1 are easily verified because

$$r_H(t) = 1 - |t|^{2H}/2 + O(t^2), t \to 0 \text{ and } r_H(t) = e^{-tH}/2 + He^{-t(1-H)} + O(e^{-t(2-H)}), t \to \infty.$$

In addition, the spectral function of $\tilde{B}_H(\tau)$ is strictly positive:

$$f_H(\lambda) = c_H \cosh(\pi\lambda)|\Gamma(-H+i\lambda)|^2 > 0$$

where $\Gamma(\cdot)$ is the Gamma-function. Therefore, $\ln f_H(\lambda) \in L_{1,loc}$; further, $f_H(\lambda) \leq f_H(0)$ because $r_H(t) > 0$. Since the Lamperti transform of $t^{H-H_i}B_{H_i}^{(i)}(t)$ is equal to $\tilde{B}_{H_i}(\tau)$, the (b,c)-conditions of Theorem 1 follow from the asymptotic behavior of $r_{H_i}(\tau)$ at zero and at infinity. Finally, the condition (1.12) of Theorem 2 means nothing else but $r_H(t) \in L_1$ and $\int r_H(t)dt > 0$.

## 2. Proof of Theorem 1.

### 2.1. The lower bound of $p_{T,n}(0)$, $1 \ll \ln n \leq CT$.

We have to estimate $p_{T,n}(0) = P\{X^{(i)}(t) - X^{(0)}(t) \leq 0, 0 \leq t \leq T, 1 \leq i \leq n\}$ for independent stationary processes $\{X^{(i)}(t), i \geq 0\}$.

Let $M_T^{(i)} = \sup\{X^{(i)}(t), 0 \leq t \leq T\}$. Since the processes $\{X^{(i)}(t), i = 0,1,\ldots,n\}$ are independent and $\{X^{(0)}(t)\} \stackrel{d}{=} \{-X^{(0)}(t)\}$, we have for any $a > 0$

$$p_{T,n}(0) \geq P\{X^{(0)}(t) \leq -a, X^{(i)}(t) \leq a, 0 \leq t \leq T, 1 \leq i \leq n\}$$

$$= P(M_T^{(0)} \leq -a)\prod_1^n P\{M_T^{(i)} \leq a\}. \tag{2.1}$$

***Estimate of*** $P(M_T^{(i)} \leq a)$.

Obviously, $P(M_T^{(i)} \leq a) \geq P(|X^{(i)}(t)| \leq a, t \in (0,T))$. For any continuous centered Gaussian process, the Gaussian correlation inequality [12, 21] implies

$$P(|X(t)| \leq a, 0 \leq t \leq T)) \geq \prod_i P\{|X(t)| \leq a, t \in \Delta_i\}, \cup \Delta_i = [0, T].$$

Assuming $\Delta_i = [i, i+1] \cdot \rho$, $T = \rho \cdot m$, and the stationarity of $X(t)$, one has

$$P(|X(t)| \leq a, 0 \leq t \leq T)) \geq [P(|X(t)| \leq a, 0 \leq t \leq \rho)]^m.$$

Recall the concentration principle for the maximum of a continuous centered Gaussian process $X(t)$, [17]. Suppose $\mu_\rho$ is the median of the distribution of $M_\rho = \max(X(t), 0 \leq t \leq \rho)$, $\sigma^2(\rho) = \max_{0 \leq t \leq \rho} E[X(t)]^2$, and $\Phi(x)$ is the standard Gaussian distribution, then for any $\tau > 0$

$$P(M_\rho \leq \mu_\rho + \tau) \geq \Phi(\tau/\sigma(\rho)). \tag{2.2}$$

Given the fact that $m_\rho = \inf(X(t), 0 \leq t \leq \rho) \stackrel{d}{=} -M_\rho$ and applying (2.2), we have

$$P(|X(t)| \leq a, 0 \leq t \leq \rho) \geq P(M_\rho \leq a) - P(m_\rho \leq -a) = P(M_\rho \leq a) - P(M_\rho \geq a)$$

$$\geq 1 - 2\Psi((a - \mu_\rho)/\sigma(\rho)), \quad a > \mu_\rho,$$

where $\Psi(x) = 1 - \Phi(x) = \Phi(-x)$.

In the stationary case, $\sigma^2(\rho) = r(0) = 1$. Since $\Psi(x) \leq 0.5 \exp(-x^2/2)$, $x > 0$, we have

$$P(|X^{(i)}(t)| \leq a, 0 \leq t \leq \rho) \geq 1 - \exp(-(a - |\mu_\rho^{(i)}|)^2/2), \quad a \geq |\mu_\rho^{(i)}| \tag{2.3}$$

Since the set of correlation functions $R$ is finite,

$$\mu_\rho = \sup_i |\mu_\rho^{(i)}| < \infty.$$

Therefore, by setting

$$a = \sqrt{2 \ln n} + \mu_\rho, \quad \rho > n/(n-1), \tag{2.4}$$

we have

$$\prod_1^n P\{M_T^{(i)} < a\} \geq [1 - \exp(-(a - \mu_\rho)^2/2)]^{nm} \geq (1 - 1/n)^{nT/\rho}$$

$$\geq (1 - 1/n)^{(n-1)T} \geq e^{-T}. \tag{2.5}$$

***Estimate of*** $P(M_T^{(0)} \leq -a)$.

Since $X^{(0)}(t)$ has a continuous spectral function $f_0(\lambda)$,

$$f_0^\bullet(\delta) = \min\{f_0(\lambda), 0 \leq \lambda \leq \delta\} \to f_0(0) = (2\pi)^{-1} \int r_0(t) dt > 0, \delta \to 0.$$

But then for small $\delta$, we can decompose the spectral function into two nonnegative terms: $f_0(\lambda) = f_\delta(\lambda) + f_c(\lambda)$ where $f_\delta(\lambda) = f_0^\bullet(\delta) 1_{|\lambda| < \delta}$. Hence, we can consider the process $X^{(0)}(t)$ as sum $X^{(0)}(t) = \xi_\delta(t) + \xi_c(t)$ of two independent stationary processes with spectral functions $f_\delta(\lambda)$ and $f_c(\lambda)$ respectively.

Given the independence of the decomposition components, we have

$$P(M_T^{(0)} \leq -a) = P\{\xi_\delta(t) + \xi_c(t) \leq -a, 0 \leq t \leq T\}$$

$$\geq P(M_T^{(\delta)} \leq -a-1) P(M_T^{(c)} \leq 1), \qquad (2.6)$$

where $M_T^{(\alpha)} = \sup\{\xi_\alpha(t), 0 \leq t \leq T\}, \alpha \in \{\delta, c\}$.

***Estimate of*** $P(M_T^{(c)} \leq 1)$.

By assumption (a),

$$E(X^{(0)}(t_1) - X^{(0)}(t_2))^2 = 2(1 - r_0(|t_1 - t_2|) \leq c|t_1 - t_2|^{2H_0}.$$

Since $E(\Delta X^{(0)})^2 = E(\Delta \xi_\delta)^2 + E(\Delta \xi_c)^2$, we have $E(\Delta \xi_c)^2 \leq c|\Delta|^{2H_0}$. Under this condition, the Talagrand's theorem [22] guarantees that

$$P\{|\xi_c(t)| \leq C; 0 \leq t \leq T\} \geq \exp(-KT/C^{1/H_0}). \qquad (2.7)$$

with some $K > 0$. In our case, $C = 1$.

***Estimate of*** $P(M_T^{(\delta)} \leq -a-1)$.

The process $\xi_\delta(\tau)$ has finite spectrum and is not correlated at the points $\tau = \pi k / \delta$. Therefore it admits the Kotelnikov-Shannon representation in terms of discrete white noise $\{\eta_n\}$:

$$\xi_\delta(\tau) = \sigma_\delta \sum_{n \in Z} \eta_n \frac{\sin(\delta\tau - \pi n)}{\delta\tau - \pi n} = \sigma_\delta S(\delta\tau / \pi), \qquad (2.8)$$

where $\sigma_\delta = \sqrt{2\delta f_0^\bullet(\delta)}$.

In terms of the random function $S(t)$ the probability under consideration is

$$P(M_T^{(\delta)} \leq -a-1) = P\{S(t) + a_\delta < 0, |t| \leq T_\delta / 2\} := Q_{T_\delta}, \qquad (2.9)$$

where $a_\delta = (a+1)/\sigma_\delta$ and $T_\delta = T\delta/\pi$.

Following [2], consider for odd $N$

$$g_N(t) = \sum_{|n| < N/2} \frac{\sin \pi(t-n)}{\pi(t-n)},$$

$$S_N(t|a) = \sum_{|n| < N/2} (\eta_n + a) \frac{\sin \pi(t-n)}{\pi(t-n)}, \qquad R_N(t) = \sum_{|n| > N/2} \eta_n \frac{\sin \pi(t-n)}{\pi(t-n)}.$$

Then

$$S(t) + a_\delta = S_N(t|a_\delta + 1) - 1 - (a_\delta + 1)(g_N(t) - 1) + R_N(t). \qquad (2.10)$$

**Lemma 4**. There exists a constant $c > 0$ such that

$$|g_N(t) - 1| < c/(N - T_\delta) + e^{-\pi N/2} := \kappa(\delta, N) \qquad (2.11)$$

for $|t| \leq T_\delta / 2 < N/2$.

**Proof of Lemma 4.** Following [2], we represent $g_N(t)$ by the contour integral

$$g_N(t) - 1 = \frac{\sin \pi t}{2\pi i} \oint_K \frac{dz}{(t-z)\sin \pi z},$$

where $K$ is the boundary of the square of size $N$, centered at the origin. On the sides of $K$ that are parallel to the $x$ axis one has

$$|t - z| \geq N/2, \quad |\sin \pi z|^2 = \cosh \pi N.$$

The analogous estimates for the sides of $K$ that are parallel to the $y$ axis are

$$|t - z| \geq (N - T_\delta)/2, \quad |\sin \pi z|^2 = \cosh(2\pi \operatorname{Im} z).$$

Substituting these estimates into the integral, we obtain (2.11).

□

Let us return to (2.10). Taking into account the independence of the terms in (2.10) and the inequality (2.11), we can estimate (2.9) as follows:

$$Q_{T_\delta} \geq P\{S_N(t|a_\delta + 1) < 1/2 - (a_\delta + 1)\kappa(\delta, N), |t| \leq T_\delta/2\}$$
$$\times P\{R_N(t) < 1/2, |t| \leq T_\delta/2\} := Q_{T_\delta}^{(1)} \times Q_{T_\delta}^{(2)}. \tag{2.12}$$

Suppose that

$$(a_\delta + 1)\kappa(\delta, N) < 1/4. \tag{2.13}$$

Then

$$Q_{T_\delta}^{(1)} \geq P\{|S_N(t|a_\delta + 1)| < 1/4, |t| \leq T_\delta/2\}.$$

The function $S_N(t|a_\delta + 1)$ is obtained from $S_N(t|0)$ by shifting the i.i.d. random variables $\{\eta_n, |n| < N/2\}$ by the constant $(a_\delta + 1)$. Therefore we can continue

$$= E 1_{\Theta_N} e^\varsigma e^{-N(a_\delta + 1)^2/2},$$

where

$$\Theta_N = \{|S_N(t|0)| < 1/4, |t| \leq T_\delta/2\} \text{ and } \varsigma = (a_\delta + 1)\sum_{|n| < N/2} \eta_n.$$

Since $\eta_n = S_N(n|0)$, the event $\Theta_N$ entails $\eta_n \geq -1/4, |n| \leq T_\delta/2$.

Therefore

$$\varsigma \geq -(a_\delta + 1)(T_\delta + 1)/4 + \varsigma_\Delta, \quad \varsigma_\Delta = (a_\delta + 1)\sum_{T_\delta/2 < |n| < N/2} \eta_n.$$

By the Cauchy-Schwartz inequality

$$(E 1_{\Theta_N} e^{\varsigma_\Delta/2} \cdot e^{-\varsigma_\Delta/2})^2 \leq E 1_{\Theta_N} e^{\varsigma_\Delta} E e^{-\varsigma_\Delta} = E 1_{\Theta_N} e^{\varsigma_\Delta} \exp((a_\delta + 1)^2(N - T_\delta)/2).$$

Hence,

$$Q_{T_\delta}^{(1)} \geq E 1_{\Theta_N} e^{\varsigma_\Delta} e^{-(a_\delta + 1)(T_\delta + 1)/4} e^{-N(a_\delta + 1)^2/2} \geq [P\{|S_N(t|0)| < 1/4, |t| \leq T_\delta/2\}]^2$$

$$\times \exp[-(2N - T_\delta)(a_\delta + 1)^2 / 2 - (a_\delta + 1)(T_\delta + 1)/4] \ . \tag{2.14}$$

Since $E|\Delta S_N(t|0)|^2 \leq E|\Delta S(t)|^2 < c|\Delta t|^2$, the Talagrand's theorem [22] implies

$$P\{|S_N(t|0)| < 1/4, |t| \leq T_\delta/2\} \geq e^{-K_1 T_\delta} = e^{-K_1 T_\delta} \tag{2.15}$$

According to [2] (section 2.2.2),

$$Q^{(2)}_{T_\delta} = P\{R_N(t) < 1/2, |t| \leq T_\delta/2\} \geq 1 - c/\sqrt{N - T_\delta} \ . \tag{2.16}$$

Putting the above inequalities (2.5),(2.7),(2.14)-(2.16) together yields

$$p_{T,n} \geq c \exp\{-KT - 2K_1 T_\delta - [2(N - T_\delta) + T_\delta](a_\delta + 1)^2/2 - (a_\delta + 1)(T_\delta + 1)/4\} \tag{2.17}$$

for $N - T_\delta \gg 1$.

To choose the appropriate $\delta$, we recall all the constants that we used:

$$T_\delta = T\delta/\pi \ , \ a_\delta = (a+1)/\sigma_\delta \ , \ \sigma_\delta = \sqrt{2\delta f_0^\bullet(\delta)} \ , \ a = \sqrt{2 \ln n} + \mu_\rho \ ,$$

and our assumptions:

$$(a_\delta + 1)[c/(N - T_\delta) + e^{-\pi \cdot N/2}] \leq 1/4 \ , \text{ and } \ N - T_\delta \gg 1. \tag{2.18}$$

Here $T, n, \delta^{-1}$ as well as $a$ and $a_\delta$ are large numbers. It is not difficult to see that we can satisfy (2.18) by setting $N - T_\delta = Ca_\delta$ with some $C > 0$.

Then the argument of the exponent in (2.17) taken with the sign (-) is

$$TK(1 + c\delta) + O(a_\delta^3) + T_\delta a_\delta^2/2 \cdot (1 + O(1/a_\delta))$$
$$= O(T) + O((\ln n/\delta)^{3/2}) + T \ln n \cdot (2\pi f_0^\bullet(\delta))^{-1}(1 + O(\sqrt{\delta/\ln n})) \ . \tag{2.19}$$

Putting $\delta^6 \ln n = 1$, we have

$$p_{T,n}(0) \geq \exp(-cT - c_1(\ln n)^{7/4} - T \ln n/(2\pi f_0^\bullet(\delta))(1 + c_2/(\ln n)^{2/3}) \ .$$

Assume that $\ln n \leq CT$, then

$$p_{T,n}(0) \geq \exp(-(T \ln n)(1 + c(T^{-1/4} + (\ln n)^{-2/3})/(2\pi f_0^\bullet(\delta))) \ , \tag{2.20}$$

where $f_0^\bullet(\delta) - f_0(0) = o(1), \delta = (\ln n)^{-1/6} \to 0$. □

**2.2. The upper bound of** $p_{T,n}(0)$, $c \ln T < \ln n$, $c > 1$.

We start from two lemmas:

**Lemma 5,[16].** Let $\{\xi_i, i = 1,..., n\}$ be a centered Gaussian stationary sequence with a correlation function $r(i)$ such that $\max_{i \geq 1} |r(i)| = \delta < r(0) = 1$. Then for any subsequence $\{\xi_{k(i)}, i = 1,..., \nu \leq n\}$

$$\left| P\{\xi_{k(i)} \leq a, i = 1,..., \nu\} - \Phi(a)^\nu \right| \leq K(1 - \delta^2)^{-1/2} \nu \sum_{i=1}^n |r(i)| \exp(-a^2/(1 + \delta)) \ . \tag{2.21}$$

**Lemma 6**. Let $X(t)$ be a centered Gaussian stationary process with a continuous spectral density $f_0(\lambda)$. Assume that

a) $f_0(\lambda) \leq f_0(0) < \infty$, b) $f_0(\lambda) \leq \varphi(\lambda) \in L_1$, where $\varphi(\lambda)$ is a monotone function,

c) $\ln f_0(\lambda) \in L_{1,loc}$.

Then for any real numbers $\{u_i\}$,

$$P\{X(i\theta) > u_i, i = 1,..., m\} \leq P\{\sigma_\theta \eta_i > u_i, i = 1,..., m\} \exp(A_\theta m / 2), \qquad (2.22)$$

where $\{\eta_i, i = 1,..., m\}$ are i.i.d. standard Gaussian variables,

$$A_\theta = \sup_{\Lambda \leq \pi/\theta} < \ln[2\pi f_0(0) / f_0(\lambda)] >_\Lambda + 1, \quad \theta \leq \theta_0, \qquad (2.23)$$

$<\psi(\lambda)>_\Lambda$ is the mean of $\psi(\lambda)$ in the interval $(0, \Lambda)$, and

$$\theta \sigma_\theta^2 = 2\pi f_0(0) + o(1), \quad \theta \to 0. \qquad (2.24)$$

**Proof of Lemma 6.** The sequence $\{X(k\theta)\}$ is stationary and has spectral function

$$f_\theta(\lambda) = \theta^{-1} f_0(\lambda/\theta) + 2\sum_{k \geq 1} \theta^{-1} f_0((\lambda + 2\pi k)/\theta), \quad |\lambda| \leq \pi. \qquad (2.25)$$

By assumptions (a,b), $f_0(\lambda/\theta) \leq f_0(0)$ and $f_0(\lambda)$ is majorized by a monotone function $\varphi(\lambda) \in L_1$. Therefore

$$\sum_{k \geq 1} f_0((\lambda + 2\pi k)/\theta) \leq \varphi(\pi/\theta) + \sum_{k \geq 1} \varphi((2k+1)\pi/\theta).$$

$$\leq 2\pi^{-1}\theta \int_{\pi/(2\theta)}^{\pi/\theta} \varphi(\lambda) d\lambda + (2\pi)^{-1}\theta \int_{\pi/\theta} \varphi(\lambda) d\lambda \leq 2\pi^{-1}\theta \int_{\pi/(2\theta)} \varphi(\lambda) d\lambda$$

Therefore, putting $\sigma_\theta^2 = \sup_{0 \leq \lambda \leq \pi} 2\pi f_\theta(\lambda)$, we get

$$\theta \sigma_\theta^2 = \sup_{0 \leq \lambda \leq \pi}[2\pi f_0(\lambda/\theta) + o(\theta)] = 2\pi f_0(0) + o(1), \quad \theta \to 0. \qquad (2.26)$$

Consider the $m \times m$ matrix $R_m = [r(i\theta - j\theta)]_{i,j=1-m}$. Because $\sigma_\theta^2 = \sup_{0 \leq \lambda \leq \pi} 2\pi f_\theta(\lambda)$, we have the following relation for the quadratic forms

$$(\mathbf{x}, R_m \mathbf{x}) = \int_{-\pi}^{\pi} \left| \sum_1^m x_k e^{ik\lambda} \right|^2 f_\theta(\lambda) d\lambda \leq \int_{-\pi}^{\pi} \left| \sum_1^m x_k e^{ik\lambda} \right|^2 \sigma_\theta^2 /(2\pi) d\lambda = (\mathbf{x}, \mathbf{x}) \sigma_\theta^2.$$

We will see later that $R_m$ is non-degenerate. Assuming the existence $R_m^{-1}$, we have

$$\exp[-(\mathbf{x}, R_m^{-1}\mathbf{x})/2] \leq \exp[-(\mathbf{x}, \mathbf{x})/2\sigma_\theta^2].$$

The last one means that for any $\{u_i\}$ we have

$$P\{X(i\theta) > u_i, i = 1,..., m\} \leq P(\sigma_\theta \eta_i > u_i, i = 1,..., m) Q_m$$

where $\{\eta_i, i = 1,..., m\}$ are i..i.d. standard Gaussian variables, $Q_m = \sqrt{\sigma_\theta^{2m} / D_m}$ and $D_m = \det R_m$.

According to the theory of the Toeplitz forms [9], $\delta_m^2 = D_m / D_{m-1}$ is the mean-square error of the $X(0)$ prediction based on $\{X(i\theta), i = 1,..., m-1\}$ data. Moreover, $\delta_m^2 = D_m / D_{m-1}$ decreases and converges to the value

$$\delta_\infty^2 = \exp\{(2\pi)^{-1} \int_{-\pi}^{\pi} \ln f_\theta(\lambda) d\lambda\}.$$

Since $f_\theta(\lambda) \geq \theta^{-1} f_0(\lambda/\theta)$ we have

$$D_m = \delta_m^2 \cdot \delta_{m-1}^2 \cdot ... \cdot \delta_1^2 \geq (\delta_\infty^2)^m \geq \exp\{m(\theta/\pi) \int_0^{\pi/\theta} \ln f_0(\lambda) d\lambda\} \cdot \theta^{-m}.$$

Therefore

$$Q_m = \sqrt{\sigma_\theta^{2m} / D_m} \leq (\theta \sigma_\theta^2)^{m/2} \exp\{-m/2 \cdot (\theta/\pi) \int_0^{\pi/\theta} \ln f_0(\lambda) d\lambda\}.$$

This estimate also proves the non-degeneracy of the matrix $R_m$. According to (2.26),

$$\theta \sigma_\theta^2 = 2\pi f_0(0) + o(1), \quad \theta \to 0.$$

Hence for small $\theta$ we have $Q_m \leq \exp(A_\theta m / 2)$ with

$$A_\theta = \exp\{(\theta/\pi) \int_0^{\pi/\theta} \ln(2\pi f_0(0) / f_0(\lambda)) d\lambda + 1.$$

Obviously, $A_\theta$ can be replaced by a non-decreasing function of $\theta^{-1}$ represented by (2.23).

□

**Proof of the upper bound.** We discretize the time in increments of $\theta$ so that $T = m\theta$ and $m$ is integer. Using notation $M_n(t) = \max[X^{(i)}(t), i = 1,..., n]$, we have

$$p_{T,n}(0) = P\{M_n(t) < X^{(0)}(t), 0 \leq t \leq T\} \leq P\{M_n(k\theta) < X^{(0)}(k\theta), k = 0,..., m\}$$

$$= EP\{M_n(k\theta) < X^{(0)}(k\theta), k = 0,..., m | M_n(\cdot)\}.$$

Lemma 5 allows replacing the process $X^{(0)}$ with a sequence $\{\sigma_\theta \eta_i, i = 1,..., m\}$ of i.i.d. Gaussian variables, where $\theta \sigma_\theta^2 = 2\pi f_0(0) + o(1), \theta \to 0$.

For this reason,

$$p_{T,n}(0) \leq P\{M_n(k\theta) < \sigma_\theta \eta_k, k = 0,..., m\} \exp(A_\theta m / 2) := J_{n,m} \exp(A_\theta m / 2), \quad (2.27)$$

where $A_\theta$ is given by (2.23).

Let $\nu_a = \#\{\sigma_\theta \eta_k > a_n, k = 1,..., m\}$, $\bar{\nu}_a = m - \nu_a$, and $\bar{\theta} = 1 - \theta$. Then

$$J_{n,m} \leq P(\nu_a > \bar{\theta} m) + P\{M_n(k\theta) < \sigma_\theta \eta_k, k = 0,..., m; \bar{\nu}_a > \theta m\} := J_{nm}^{(1)} + J_{nm}^{(2)}.$$

Since $\nu_a$ is a binomial random variable with parameters $m$ and $p = \Psi(a_n / \sigma_\theta)$, we have

$$J_{nm}^{(1)} \leq \sum_{k > \bar{\theta} m} C_m^k p^k (1-p)^{m-k} \leq C_m^{[m/2]} p^{[\bar{\theta} m]} (1-p)^{-1} \leq c 2^m [\Psi(a_n / \sigma_\theta)]^{\bar{\theta} m},$$

$$J_{nm}^{(1)} \geq p^{[\bar{\theta} m]} (1-p)^m > 2^{-m} [\Psi(a_n / \sigma_\theta)]^{\bar{\theta} m}, \text{ (if } p \leq 1/2 \text{).}$$

Note, that $\theta\sigma_\theta^2 = d_0 + o(1), \theta \to 0$, $d_0 = 2\pi f_0(0)$ ; $T = \theta \cdot m$ , and

$$c < \Psi(u)/[u^{-1}e^{-u^2/2}] < C, u > u_0 > 0 \ . \qquad (2.28)$$

Here and below $\theta$ and the other parameters depend on $n$ and therefore the notation $o(1)$ means that $o(1) \to 0$ as $n \to \infty$ .

For any $a_n^2 = 2\ln n(1 + o(1))$ we have

$$J_{nm}^{(1)} \le c \exp\{-\overline{\theta} T \ln n/(d_0 + o(1)) + T\theta^{-1} \ln 2\} ,$$

$$J_{nm}^{(1)} \ge c \exp\{-\overline{\theta} T \ln n/(d_0 + o(1)) - T\theta^{-1} \ln(c_1 \ln n)\} .$$

Hence

$$J_{nm}^{(1)} = \exp\{-T\ln n/(d_0 + o(1))\} \text{ if } 1/\theta = o(\ln n/\ln\ln n) , \theta = o(1) . \qquad (2.29)$$

Moreover

$$\ln[J_{nm}^{(1)} \exp(A_\theta T/2\theta)] = -T\ln n/d_0(1 + o(1)) \text{ if } A_\theta/\theta = o(\ln n) . \qquad (2.30)$$

***Estimate of*** $J_{nm}^{(2)}$.

Let $I_a = \{k_i, i = 1,..., \overline{v}_a : \eta_{k_i}\sigma_\theta \le a_n\}$ . Then

$$J_{nm}^{(2)} = P\{M_n(k\theta) < \sigma_\theta \eta_k, k = 0,..., m; \overline{v}_a \ge \theta m\}$$

$$\le P\{X^{(p)}(k_i\theta) < a_n, p = 1,..., n; k_i \in I_a, \overline{v}_a \ge \theta m\}$$

$$= E\prod_{p=1}^n P\{X^{(p)}(k_i\theta) < a_n | \{k_i\} = I_a, \overline{v}_a \ge \theta \cdot m\} .$$

Applying Lemma 5, we can continue

$$\le E\Pi_{i=1}^n[\Phi^{\overline{v}_a}(a_n) + c(1 - \delta_{\theta,i}^2)^{-1/2}\overline{v}_a\theta^{-1}\exp(-a_n^2/(1+\delta_{\theta,i}), i = 1,..., \overline{v}_a)|\overline{v}_a \ge \theta \cdot m = T]$$

Since the correlation functions set $R$ is finite, we have

$$\sum_{i\ge 1}|r_p(i\theta)|\theta \le \int_0^\infty \tilde{r}_p(t)dt < const \text{ and } \delta_\theta = \sup_{i\ge 1}\delta_{\theta,i} < 1 \ .$$

Therefore

$$J_{nm}^{(2)} \le [\Phi^T(a_n) + c(1 - \delta_\theta^2)^{-1/2}T\theta^{-1}\exp(-a_n^2/(1+\delta_\theta))]^n := [W_1 + W_2]^n .$$

Now we will refine $a_n$ by finding it as the root of the equation $\Psi(a_n) = n^{-1+\vartheta}, 0 < \vartheta < 1$ .

It's easy to see that

$$a_n^2 = (1 - \vartheta)2\ln n - 2\ln\ln n \cdot (1 + o(1)) ,$$

i.e. $a_n^2 = 2\ln n(1 + o(1))$ if $\vartheta = o(1)$. In this case

$$W_1^n = [\Phi^T(a_n)]^n = (1 - n^\vartheta/n)^{nT} \le \exp(-n^\vartheta T) \ . \qquad (2.31)$$

Now we estimate $[1 + W_2/W_1]^n$. First of all,

$$W_1 = (1 - n^\vartheta / n)^T \geq \exp(-T/(n^{1-\vartheta} - 1))$$

$$\geq \exp(-2n^{\vartheta+q-1}), \quad \text{if } \ln T \leq q \ln n, \quad 0 < \vartheta + q < 1, \quad n > n_0 \tag{2.32}$$

$$W_2 = c(1 - \delta_\theta^2)^{-1/2} T \theta^{-1} \exp(-a_n^2/(1 + \delta_\theta))$$

Let $\bar{\delta}_\theta = 1 - \delta_\theta$ and $\varphi(\theta) = \bar{\delta}_\theta^{1/2} \theta$, then

$$W_2 / W_1 \leq cT / \varphi(\theta) \cdot \exp[-(1-\vartheta)/(1-\bar{\delta}_\theta/2) \ln n - 2 \ln \ln n (1 + o(1)) + 2n^{\vartheta+q-1}]$$

$$\leq n^{-1} cT / \varphi(\theta) \cdot \exp(-\ln n \cdot (\bar{\delta}_\theta/2 - \vartheta)/(1 - \bar{\delta}_\theta/2) + 1), \quad n \geq n_0.$$

By setting $\vartheta = \bar{\delta}_\theta / 4$, we get

$$[1 + W_2 / W_1]^n \leq \exp(\tilde{c} T / \varphi(\theta) \cdot n^{-\bar{\delta}_\theta/4}).$$

Given (2.31), we have

$$J_{nm}^{(2)} \leq [W_1 + W_2]^n \leq \exp[-T(n^{\bar{\delta}_\theta/4} - \tilde{c} n^{-\bar{\delta}_\theta/4} \bar{\delta}_\theta^{-1/2} \theta^{-1})].$$

Under conditions (b) and finiteness of correlation functions set, we can choose such constants, particularly $h = \sup_{i \geq 1} h_i$, that for any $i \geq 1$

$$1 - |r_i(t)| \geq c|t|^{2h} \wedge \rho, \quad h > 0, 0 < \rho < 1.$$

Therefore, for small enough $\theta$, we have $\bar{\delta}_\theta \geq c\theta^{2h}$. Assuming

$$\theta^{-1} \leq \ln^\kappa n, \quad 2h\kappa < 1, \tag{2.33}$$

we will have $\bar{\delta}_\theta \geq c\theta^{2h} \geq c(\ln n)^{-2h\kappa}$. It follows that

$$\tilde{c} n^{-\bar{\delta}_\theta/4} \bar{\delta}_\theta^{-1/2} \theta^{-1} \leq c \exp(-c_1 (\ln n)^{1-2h\kappa})(\ln n)^{(1+h)\kappa} = o(1).$$

Hence $J_{nm}^{(2)} \leq \exp[-T(n^{\bar{\delta}_\theta/4} - o(1))] = o(J_{nm}^{(1)})$.

It remains only to fulfill the requirements (2.29) and (2.30).i.e., $1/\theta = o(\ln n / \ln \ln n)$, $A_\theta / \theta = o(\ln n)$. By definition (see 2.23), $A_\theta / \theta$ is unlimited strictly increasing function of $\theta^{-1}$. Therefore we can choose $\theta = \theta_n$ as the root of the equation $A_\theta / \theta = (\ln n)^\kappa$. As a result, we get $1/\theta \leq A_\theta / \theta = (\ln n)^\kappa$. The proof is complete.

## 3. Proof of Theorem 2.

The argument of the proof is closely related to our work [18], but differs in some details due to the specifics of the problem under consideration.

***The lower bound.*** Let $\{X^{(i)}(t), X^{(i)}(0) = 0, i = 0, 1, ..., n\}$ be independent centered Gaussian continuous H-ss processes; $M_n(\Delta)$ is the maximum of the process $M_n(t) = \max\{X^{(i)}(t) - X^{(0)}(t), 1 \leq i \leq n\}$ on the interval $\Delta$ and $G_n(\Delta)$ is the far right position of $M_n(\Delta)$ on $\Delta$. Assume that

$$\liminf_{T,n\to\infty} \ln P(M_n(1,T) \leq 0)/(\ln T \cdot \ln n) = \gamma_- > -\infty , \qquad (3.1)$$

where $(n,T)$ increase such that $\ln n < C \ln T$.

For any $c_{T,n} > 0$,

$$P(M_n([1,T]) \leq 0) \leq P(G_n([0,T]) \leq 1) \leq P(G_n([0,T]) \leq 1, M_n([0,1]) \leq c_{T,n})$$
$$+ P(M_n([0,1]) \geq c_{T,n}) \leq P(M_n([0,T]) \leq c_{T,n}) + P(M_n([0,1]) \geq c_{T,n}). \qquad (3.2)$$

The H-ss property of the processes in question entails

$$P(M_n([0,T]) \leq c_{T,n}) = P(M_n([0,T']) \leq 1), \qquad T' = Tc_{T,n}^{-1/H}. \qquad (3.3)$$

At the same time,

$$R_{T,n} := P(M_n([0,1]) \geq c_{T,n}) \leq \sum_{i=1}^{n} P(\max[\xi^{(i)}(t), 0 \leq t \leq 1] \geq c_{T,n}),$$

where $\xi^{(i)}(t) = X^{(i)}(t) - X^{(0)}(t)$.

Let $\mu_i$ be the median of the distribution of $M^{(i)} = \max(\xi^{(i)}(t), 0 \leq t \leq 1))$, and

$$\sigma_i^2 = \max_{0 \leq t \leq 1} E[\xi^{(i)}(t)]^2 = \max_{0 < t < 1} t^{2h} E[\xi^{(i)}(1)]^2 = E[\xi^{(i)}(1)]^2.$$

Obviously, we can assume that $\sup \mu_i \leq \mu < \infty$ and $\sup \sigma_i^2 \leq \sigma^2 < \infty$ since the set of the correlation functions of $\{X^{(i)}(t)\}$ is finite. Therefore, applying the concentration principle to $\{M^{(i)}\}$, [17], we get for large $c_{n,T}$

$$R_{n,T} \leq \sum_{1}^{n} \Psi((c_{T,n} - \mu_i)/\sigma_i) \leq n\Psi((c_{T,n} - \mu)/\sigma), \qquad (3.4)$$

Put $(c_{T,n}/\sigma)^{2.}/2 = A \ln T \cdot \ln n$, $A > |\gamma_-|$. Then for large $T, n$

$$R_{T,n} \leq n \exp((c_{T,n} - \mu)^2/2) = \exp[-A \ln T \cdot \ln n(1 + o(1))].$$

Hence, by (3.1)

$$R_{T,n} \leq o(P(M_n([1,T]) \leq 0)). \qquad (3.5)$$

Collecting (3.2), (3.3) and (3.5) together we have

$$\frac{\ln P(M_n([0,T']) \leq 1)}{\ln T' \cdot \ln n} \geq \frac{\ln P(M_n(1,T) \leq 0)(1 + o(1))}{\ln T \cdot \ln n} \times \frac{\ln T}{\ln T'}.$$

Assuming $c < \ln n < C \ln T$ and taking $T' = cT [\ln n \cdot \ln T]^{-1/2H}$ into account, we get $\ln T / \ln T' = 1 + o(1)$ and

$$\liminf_{T',n\to\infty} \frac{\ln P(M([0,T']) \leq 1)}{\ln T' \cdot \ln n} \geq \gamma_-, \quad \ln n < C' \ln T'.$$

## *The upper bound.*

Now we estimate $P(M_n([0,T]) \leq 1)$ from above using the upper bound $\gamma_+$ of $\ln P(M_n([1,T]) \leq 0)/\psi_{T,n}$, where $\psi_{T,n} = \ln T \cdot \ln n$.

Let $\mu(s) \geq 1$, $1 \leq s \leq T$ be an element of the reproducing kernel Hilbert space $H_0(T)$, associated with $X^{(0)}(t), 0 \leq t \leq T$ (see [17]). We will need the following notation:

$$X_\mu^{(0)}(t) = X^{(0)}(t) + \mu(t),$$

$$M_n(t|\mu) = \max\{X^{(i)}(t) - X_\mu^{(0)}(t), 1 \leq i \leq n\}, \text{ and } M_n(\Delta|\mu) = \max_{t \in \Delta} M_n(t|\mu).$$

Obviously, $M_n(\Delta|0) = M_n(\Delta)$. One has

$$P\{M_n([0,T]|0) \leq 1\} = P\{M_n(t|\mu) \leq 1 - \hat{\mu}(t), 0 \leq t \leq T)\}$$

$$\leq P\{M_n([1,T])|\mu) \leq 0\} := P(\Omega_\mu) \qquad (3.6)$$

The function $\mu(s)$ is the admissible shift of the Gaussian measure $P^{(0)}(d\omega)$ related to the process $X^{(0)}(t)$ on the interval $[0,T]$, [17]. Therefore

$$P(\Omega_\mu) = E 1_{\Omega_0} \pi(\omega_T),$$

where $\pi(\omega_T)$ is the Radon-Nikodym derivative of two Gaussian measures corresponding to the processes $X^{(0)} + \mu$ and $X^{(0)}$ on $[0,T]$. Using $H\ddot{o}lder$'s inequality with parameters $(p,q) = ((1-\varepsilon)^{-1}, \varepsilon^{-1})$, we get,

$$E 1_{\Omega_\mu} \pi(\omega_T) \leq [P(M_n((1,T)|0)]^{1-\varepsilon} \exp(\varepsilon \psi_{T,n}), \qquad (3.7)$$

where $\varepsilon^2 = \|\mu\|^2 / \psi_{T,n}$ and $\|\cdot\|$ is the norm in $H_0(T)$ (for more information see [18]). Suppose that $\varepsilon = o(1)$ then

$$\ln P(M_n((0,T)|\mu) \leq 0) / \psi_{n,T} \leq (1 - o(1)) \ln P(M_n((1,T)|0) / \psi_{n,T} + o(1).$$

It remains to find an element $\mu(t)$ from the Hilbert space $H_0(T)$ such that $\mu(s) \geq 1$, $1 \leq s \leq T$ and $\|\mu\|^2 = o(\ln T \ln n)$.

Let us consider the Lamperti transform, $\tilde{X}^{(0)}(\tau) = X^{(0)}(e^\tau) e^{-\tau \cdot h} := LX^{(0)}(\tau)$, of the process $X^{(0)}(t)$. Assume that $\tilde{H}_0(\tilde{T})$ is the Hilbert spaces with the reproducing kernels associated with the process $\tilde{X}^{(0)}(\tau), \tau \in (-\infty, \tilde{T} = \ln T)$. Then the mapping $\psi \to \tilde{\psi} = L\psi$ performs isometry of spaces $H_0(T)$ and $\tilde{H}_0(\tilde{T})$.

Assume that $dF_0(\lambda)$ is the spectral measure of the stationary process $\tilde{X}^{(0)}(\tau) = LX^{(0)}$. Then, as noted in [8],

$$\tilde{\psi}(\tau) = 2\int_{0+}^{1/\tilde{T}} \cos \tau\lambda \, dF_0(\lambda) / F_0((0, \tilde{T}^{-1}]) \geq 1, \tau \in [0, \tilde{T}] \text{ and } \|\tilde{\psi}\|_{\tilde{H}}^2 = 2 / F_0((0, \tilde{T}^{-1}])$$

where $\tilde{H} = \tilde{H}_0(\infty)$. Conditions (1.12) mean that the spectral measure $dF_0(\lambda)$ has continuous density $f_0(\lambda), f_0(0) > 0$. Therefore $F_0((0, \tilde{T}^{-1}]) \approx f_0(0) / \tilde{T}$ and

$$\psi(t) := L^{-1}\tilde{\psi} = \tilde{\psi}(\ln t) \cdot t^H \geq 1 \;,\; 1 \leq t \leq T \;,\; \|\psi\|_H^2 = \|\tilde{\psi}\|_{\tilde{H}}^2 \leq C\tilde{T} = C \ln T \;.$$

Assume that $\mu(t)$ is projection of $\psi(t)$ on $H_0(T) \subset H_0(\infty) := H$. Then $\mu(t) = \psi(t)$ on $(0,T)$ and

$$\|\mu\|_{H_0(T)}^2 \leq \|\psi\|_H^2 \leq C \ln T = o(\ln n \cdot \ln T).$$

The element $\mu(t)$ is the desired function from $H_0(T)$.

**Acknowledgements**. This research was supported by the Russian Science Foundation through the Research Project 17-11-01052.